# A LINEAR, DECOUPLED AND ENERGY STABLE SCHEME FOR SMECTIC–A LIQUID CRYSTAL FLOWS

XIAOFENG YANG [†⋆] AND ALEX BRYLEV [‡]

Abstract. In this paper, we consider numerical approximations for the model of smectic-A liquid crystal flows. The model equation, that is derived from the variational approach of the de Gennes free energy, is a highly nonlinear system that couples the incompressible Navier-Stokes equations, and two nonlinear coupled second-order elliptic equations. Based on some subtle explicit–implicit treatments for nonlinear terms, we develop a unconditionally energy stable, linear and decoupled time marching numerical scheme. We also rigorously prove that the proposed scheme obeys the energy dissipation law at the discrete level. Various numerical simulations are presented to demonstrate the accuracy and the stability thereafter.

## 1. Introduction

Liquid crystal (LC) is often viewed as the fourth state of the matter besides the gas, liquid and solid. It may flow like a liquid, but its molecules may be oriented in a crystal-like way. There are many different types of liquid-crystal phases, which can be distinguished by their different optical properties. Thermotropic LCs can be distinguished into two main different phases: Nematic and Smectic. In Nematic phases, the rod-shaped molecules have no positional order, but molecules self-align to have a long-range directional order with their long axes roughly parallel. Thus, the molecules are free to flow and their center of mass positions are randomly distributed as in a liquid, although they still maintain their long-range directional order. In smectic phases, which are found at lower temperatures, the well-defined layers form, that can slide over one another in a manner similar to that of soap. The smectics are thus positionally ordered along one direction inside the layer. There are many different smectic phases, all characterized by different types and degrees of positional and orientational order (cf. [3,7,11,12,17,37,39,42–44,44,47–49] and the references therein). In particular, in Smectic-A phase, molecules are oriented along the normal vector of the layers, while in Smectic-C phases they are tilted away from the normal vector of the layer.

---

[†⋆] Corresponding author, Department of Mathematics, University of South Carolina, Email: xfyang@math.sc.edu.

[‡] Department of Mathematics, University of South Carolina, Email: brylev@math.sc.edu.

The mathematical model of liquid crystals can often be derived from an energy-based variational formalism (energetic variational approach), leading to well-posed nonlinearly coupled systems that satisfy thermodynamics-consistent energy dissipation laws. This makes it possible to carry out mathematical analysis and design numerical schemes which satisfy a corresponding energy dissipation law at the discrete level. For smectic-A phase liquid crystals, in de Gennes' pioneering work [7], the phenomenonlogical free energy of smectic–A phase is presented by coupling two order parameters which characterize the average direction of molecular alignment, as well as the layer structure, respectively. In [4], the authors modified the de Gennes' model by adding a second order gradient term for the smectic order parameter to investigate the nematic to smectic–A or smectic–C phase transition, and to predict the twist grain boundary phase in chiral smectic liquid crystals. In [18], the authors used the de Gennes energy to study smectic–A liquid crystals to simulate the chevron (zigzag) pattern formed in the presence of an applied magnetic field. In [8], the authors derived the hydrodynamics coupled model for smectic–A phase by assuming that the director field is strictly equal to the gradient of the layer, thus the free energy is reduced to one order parameter.

From the numerical point of view, it is specifically desired to design numerical schemes that could preserve the thermo-dynamically consistent dissipation law (energy-stable) at the discrete level, since the preservation of such laws is critical for numerical methods to capture the correct long time dynamics. The noncompliance of energy dissipation laws may lead to spurious numerical solutions if the grid and time step sizes are not carefully controlled. To the best of the author's knowledge, although a variety of the smectic liquid crystal models had been developed, we notice that the successful attempts in designing efficient energy stable schemes are very scarce due to the complex nonlinearities. For instances, in [14], the authors present a temporal second order numerical scheme to solve the model of [8]. The scheme is energy stable, however, it is nonlinear thus the implementation is complicated and the computational cost might be high. In [18], the authors developed a temporal first order scheme. However, it does not follow the energy dissipation law even though the schemes are linear and decoupled.

Therefore, the main purpose of this paper is to construct the efficient schemes to solve the de Gennes type smectic–A liquid crystal model (cf. [7,18]). We first couple the hydrodynamics to the original de Gennes free energy and derive the whole model based on the variational approach and the Ficks' law. To solve the model, the main difficulties roughly include (i) the coupling between the velocity and director field/layer function through the convection terms and nonlinear stresses; (ii) the coupling of the velocity and pressure through the incompressibility constraint; (iii) the nonlinear coupling between the director field and layer function. We develop a time discretization scheme which (a) is unconditionally stable; (b) satisfies a discrete energy law; and (c) leads to linear, decoupled equations to solve at each time step. This is by no means an easy task due to many highly nonlinear terms and the couplings existed in the model.

The rest of the paper is organized as follows. In Section 2, we present the whole model and present the PDE energy law. In Section 3, we develop the numerical scheme and prove the unconditional stability. In Section 4, we present some numerical experiments to validate the proposed scheme. Finally, some concluding remarks are presented in Section 5.

## 2. The smectic-A liquid crystal fluid flow model and its energy law

The de Gennes free energy of smectic A liquid crystal is described by a unit vector (director field) $\boldsymbol{d}$ and a complex order parameter $\psi$, to represent the average direction of molecular alignment and the layer structure, respectively. The smectic order parameter is written as

$$\psi(\boldsymbol{x}) = \rho(\boldsymbol{x})e^{iq\omega(\boldsymbol{x})}, \tag{2.1}$$

where $\omega(\boldsymbol{x})$ is the order parameter to describe the layer structure so that $\nabla\omega$ is perpendicular to the layer. The smectic layer density $\rho(\boldsymbol{x})$ is the mass density of the layers. The de Gennes free energy reads as follows,

$$E(\psi, \boldsymbol{d}) = \int_\Omega \Big( C|\nabla\psi - iq\boldsymbol{d}\psi|^2 + K|\nabla \boldsymbol{d}|^2 + \frac{g}{2}(|\psi|^2 - \frac{r}{g})^2 - \chi_a H^2(\boldsymbol{d}\cdot\boldsymbol{h})^2\Big) d\boldsymbol{x}, \tag{2.2}$$

where the order parameters $C, K, g, r$ are all fixed positive constants, $\boldsymbol{h}$ is the unit vector representing the direction of magnetic field and $H^2$ is the strength of the applied field. $\Omega = (-L, L)^2 \times (-d, d)$.

Now we consider the simple case by assuming the density $\rho(x) = r/g$ [18], then the energy becomes

$$E(\psi, \boldsymbol{d}) = \int_\Omega \Big( Cq^2|\nabla\omega - \boldsymbol{d}|^2 + K|\nabla \boldsymbol{d}|^2 - \chi_a H^2(\boldsymbol{d}\cdot\boldsymbol{h})^2\Big) d\boldsymbol{x}. \tag{2.3}$$

Let $\phi(\boldsymbol{x}) = \frac{\omega(\boldsymbol{x})}{d}$, thus the normalized energy becomes

$$E(\phi, \boldsymbol{d}) = \lambda \int_{\tilde\Omega} \Big( \frac{1}{\eta}\frac{|\nabla\phi - \boldsymbol{d}|^2}{2} + \eta\frac{|\nabla \boldsymbol{d}|^2}{2} - \frac{\tau}{2}(\boldsymbol{d}\cdot\boldsymbol{h})^2\Big) d\boldsymbol{x}, \tag{2.4}$$

where

$$\begin{aligned} &\tilde{\boldsymbol{x}} = \frac{\boldsymbol{x}}{d},\ \tilde\Omega = (0, 2\ell)^2 \times (0, 2),\ \ell = \frac{L}{d}, \\ &\eta = \frac{\gamma}{d}, \gamma = \sqrt{\frac{K}{Cq^2}},\ \lambda = \frac{2dK}{\eta},\ \tau = \frac{\chi_a H^2 d^2 \eta}{K}. \end{aligned} \tag{2.5}$$

The dimensionless parameter $\eta$ is in fact the ratio of the layer thickness to the sample thickness and thus $\eta \ll 1$.

To release the unit vector constraint of $|\boldsymbol{d}| = 1$, a nonlinear potential $G(\boldsymbol{d}) = \frac{1}{4\epsilon^2}(|\boldsymbol{d}|^2 - 1)^2$ which is a Ginzburg-Landau type penalty term, is added into the free energy to approximate the unit length constraint of $\boldsymbol{d}$ [20, 21] , where $\epsilon \ll 1$ is a penalization parameter. Thus the modified total free energy with the hydrodynamics reads as follows

$$(2.6)\quad E_{tot}(\boldsymbol{u}, \phi, \boldsymbol{d}) = \int_{\tilde{\Omega}} \frac{1}{2}|\boldsymbol{u}|^2 d\boldsymbol{x} + \lambda \int_{\tilde{\Omega}} \Big(\frac{1}{\eta}\frac{|\nabla\phi - \boldsymbol{d}|^2}{2} + \eta\frac{|\nabla \boldsymbol{d}|^2}{2} + \mathbf{G}(\boldsymbol{d}) - \frac{\tau}{2}(\boldsymbol{d}\cdot\boldsymbol{h})^2\Big) d\boldsymbol{x},$$

where $\boldsymbol{u}$ is the fluid velocity field.

Assuming a generalized Fick's law that the mass flux be proportional to the gradient of the chemical potential [2, 10, 22, 24], we can derive the following system:

$$(2.7)\quad \phi_t + \boldsymbol{u}\cdot\nabla\phi = -M_1\mu_\phi,\ \mu_\phi = \frac{\delta E}{\delta \phi},$$

$$(2.8)\quad \boldsymbol{d}_t + \boldsymbol{u}\cdot\nabla\boldsymbol{d} = -M_2\boldsymbol{\mu}_d,\ \boldsymbol{\mu}_d = \frac{\delta E}{\delta \boldsymbol{d}},$$

$$(2.9)\quad \boldsymbol{u}_t + \boldsymbol{u}\cdot\nabla\boldsymbol{u} - \nabla\cdot\sigma_d + \nabla p - \mu_\phi\nabla\phi - \mu_d\nabla\boldsymbol{d} = 0,$$

$$(2.10)\quad \nabla\cdot\boldsymbol{u} = 0,$$

where $p$ is the pressure, $\sigma_d$ is the Caughy stress tensor, $\nu$ is the viscosity, $M_1, M_2$ are the relaxation order parameters. The variational derivatives $\mu_\phi$ and $\boldsymbol{\mu}_d$ are

$$(2.11)\quad \mu_\phi = \lambda\frac{1}{\eta}(-\Delta\phi + \nabla\cdot\boldsymbol{d}),$$

$$(2.12)\quad \boldsymbol{\mu}_d = \lambda\Big( - \eta\Delta\boldsymbol{d} + \mathbf{g}(\boldsymbol{d}) + \frac{1}{\eta}(-\nabla\phi + \boldsymbol{d}) - \tau(\boldsymbol{d}\cdot\boldsymbol{h})\boldsymbol{h}\Big),$$

where $\mathbf{g}(\boldsymbol{d}) = \frac{1}{\epsilon^2}\boldsymbol{d}(|\boldsymbol{d}|^2 - 1)$.

Following the work in [8], the Cauchy stress tensor $\sigma_d$ reads as follows,

$$(2.13)\quad \sigma_d = \mu_1(\boldsymbol{d}^T D(\boldsymbol{u})\boldsymbol{d})\boldsymbol{d}\otimes\boldsymbol{d} + \mu_4 D(\boldsymbol{u}) + \mu_5(D(\boldsymbol{u})\boldsymbol{d}\otimes\boldsymbol{d} + \boldsymbol{d}\otimes D(\boldsymbol{u})\boldsymbol{d}),$$

where $\mu_i > 0$ and $D(\boldsymbol{u}) = \frac{1}{2}((\nabla\boldsymbol{u}) + (\nabla\boldsymbol{u})^T)$ is the strain tensor. In this paper, we assume $\mu_1, \mu_5$ are negligible comparing to $\mu_4$, thus the stress tensor is simplifed to

$$(2.14)\quad \sigma_d = \mu_4 D(\boldsymbol{u}).$$

For simplicty, we assume the boundary conditions as follows.

$$(2.15)\quad \boldsymbol{u}|_{\partial\Omega} = 0, \partial_{\boldsymbol{n}}\phi|_{\partial\Omega} = 0, \partial_{\boldsymbol{n}}\boldsymbol{d}|_{\partial\Omega} = 0,$$

where $\boldsymbol{n}$ is the outward normal of the boudary.

To obtain the dissipation law of the system (2.7)-(2.9), we take the $L^2$ inner product of (2.7) with $\frac{\delta E}{\delta\phi}$, (2.8) with $\frac{\delta E}{\delta \boldsymbol{d}}$, and (2.9) with $\boldsymbol{u}$, perform the integration by parts, and add all equalities together, we obtain

$$(2.16)\quad \frac{d}{dt}E_{tot}(\boldsymbol{u}, \phi, \boldsymbol{d}) = -\int_\Omega \Big(\mu_4|D(\boldsymbol{u})|^2 + M_1|\frac{\delta E}{\delta\phi}|^2 + M_2|\frac{\delta E}{\delta\boldsymbol{d}}|^2\Big) d\boldsymbol{x} \le 0.$$

## 3. Numerical scheme

The emphasis of our algorithm development is placed on designing numerical schemes that are not only easy-to-implement, but also satisfy a discrete energy dissipation law. We will design schemes that in particular can overcome the following difficulties, namely,

- the coupling of the velocity and pressure through the incompressibility condition;
- the stiffness in the director equation associated with the penalty parameter $\epsilon$;
- the nonlinear couplings among the fluid equation, the layer equation and the director equation.

We construct an energy stable scheme based on a stabilization approach [30, 31]. To this end, we shall assume that $\mathbf{G}(\boldsymbol{d})$ satisfies the following conditions, i.e.,

$$|H_{\mathbf{G}}(\boldsymbol{x})| \leq L, \forall \boldsymbol{x} \in \mathbb{R}^3. \tag{3.1}$$

where $(H_{\mathbf{G}}(\boldsymbol{x}))_{i,j} = \frac{\partial^2 \mathbf{G}}{\partial x_i \partial x_j}, i, j = 1, 2, 3$ is the Hessian matrix of $\mathbf{G}(\boldsymbol{x})$. One immediately notes that this condition is not satisfied by this double-well potential $\mathbf{G}(\boldsymbol{x})$. However, it is a common practice that one truncates this fourth order polynomial $\mathbf{G}$ to quadratic growth outside of an interval $[-M, M]$ without affecting the solution if the maximum norm of the initial condition $|\boldsymbol{d}_0|$ is bounded by $M$. Therefore, one can (cf. [6, 19, 30]) consider the truncated double-well potential $\tilde{\mathbf{G}}(\boldsymbol{d})$ by modifying this function outside a ball in $\{\boldsymbol{x} : |\boldsymbol{d}(\boldsymbol{x})| \leq 1\} \in \mathbb{R}^3$ of radius 1 as follows.

$$\tilde{\mathbf{G}}(\boldsymbol{d}) = \begin{cases} \dfrac{1}{4\epsilon^2}(|\boldsymbol{d}|^2 - 1)^2, & |\boldsymbol{d}| \leq 1, \\ \dfrac{1}{4\epsilon^2}(|\boldsymbol{d}| - 1)^2, & |\boldsymbol{d}| > 1. \end{cases} \tag{3.2}$$

Hence, there exists a postive constant $L_G$ such that

$$\max_{\boldsymbol{x} \in \mathbb{R}^3} |H_{\tilde{\mathbf{G}}}(\boldsymbol{x})| \leq L. \tag{3.3}$$

For convenience, we consider the problem formulated with the substitute $\tilde{\mathbf{G}}$, but omit the $\tilde{}$ in the notation.

When deriving the energy law (2.16), we notice that the nonlinear terms in $\frac{\delta \mathscr{E}_{tot}}{\delta \phi}$ and $\frac{\delta \mathscr{E}_{tot}}{\delta \boldsymbol{d}}$ involve second order derivatives, and it is not convenient to use them as test functions in numerical approximations, making it difficult to prove the discrete energy dissipation law. To overcome this difficulty, we first reformulate the system (2.7)-(2.10) in an alternative form which is convenient for numerical approximations. The system reads as follows,

$$\phi_t + \boldsymbol{u} \cdot \nabla \phi = -M_1 \mu_\phi, \tag{3.4}$$

$$\boldsymbol{d}_t + \boldsymbol{u} \cdot \nabla \boldsymbol{d} = -M_2 \boldsymbol{\mu}_d, \tag{3.5}$$

$$\boldsymbol{u}_t + \boldsymbol{u} \cdot \nabla \boldsymbol{u} - \nabla \cdot \sigma_d + \nabla p + \frac{\dot{\phi}}{M_1} \nabla \phi + \frac{\dot{\boldsymbol{d}}}{M_2} \nabla \boldsymbol{d} = 0, \tag{3.6}$$

$$\nabla \cdot \boldsymbol{u} = 0, \tag{3.7}$$

where $\dot{\phi} = \phi_t + \boldsymbol{u}\cdot\nabla\phi$ and $\dot{\boldsymbol{d}} = \boldsymbol{d}_t + \boldsymbol{u}\cdot\nabla\boldsymbol{d}$. To obtain the dissipation law of the system (3.4)-(3.7), we take the $L^2$ inner product of (3.4) with $\phi_t$, (3.5) with $\boldsymbol{d}_t$, and (3.6) with $\boldsymbol{u}$, perform the integration by parts, and add all equalities together. We obtain

$$\frac{d}{dt}E_{tot} = -\int_\Omega\Big(\mu_4|\nabla\boldsymbol{u}|^2 + \frac{1}{M_1}|\dot{\phi}|^2 + \frac{1}{M_2}|\dot{\boldsymbol{d}}|^2\Big)d\boldsymbol{x} \le 0. \tag{3.8}$$

We now fix some notations. For scalar function $u, v$ and vector function $\boldsymbol{u} = (u_1, u_2, u_3)$ and $\boldsymbol{v} = (v_1, v_2, v_3)$, we denote the $L^2$ inner product as follows.

$$(u,v) = \int_\Omega uv dx,\ \|u\|^2 = (u,u),\ (\boldsymbol{u},\boldsymbol{v}) = \int_\Omega \boldsymbol{u}\boldsymbol{v}^T d\boldsymbol{x},\ \|\boldsymbol{u}\|^2 = (\boldsymbol{u},\boldsymbol{u}). \tag{3.9}$$

Now, we are ready to present our energy stable scheme that reads as follows.

Given the initial conditions $\phi^0$, $\boldsymbol{d}^0$, $\boldsymbol{u}^0$ and $p^0 = 0$, having computed $\phi^n$, $\boldsymbol{d}^n$, $\boldsymbol{u}^n$ and $p^n$ for $n > 0$, we compute $(\phi^{n+1}, \boldsymbol{d}^{n+1}, \tilde{\boldsymbol{u}}^{n+1}, \boldsymbol{u}^{n+1}, p^{n+1})$ by

*Step 1:*

$$\begin{cases} \dfrac{1}{M_1}\dot{\phi}^{n+1} = \dfrac{\lambda}{\eta}(\Delta\phi^{n+1} - \nabla\cdot\boldsymbol{d}^n), \\ \dfrac{\partial\phi^{n+1}}{\partial\boldsymbol{n}}|_{\partial\Omega} = 0, \end{cases} \tag{3.10}$$

with

$$\dot{\phi}^{n+1} = \frac{\phi^{n+1}-\phi^n}{\delta t} + (\boldsymbol{u}_\star^n\cdot\nabla)\phi^n,\ \boldsymbol{u}_\star^n = \boldsymbol{u}^n - \delta t\frac{\dot{\phi}^{n+1}}{M_1}\nabla\phi^n. \tag{3.11}$$

*Step 2:*

$$\begin{cases} S(\boldsymbol{d}^{n+1} - \boldsymbol{d}^n) + \dfrac{1}{M_2}\dot{\boldsymbol{d}}^{n+1} \\ \qquad = \lambda\Big(\eta\Delta\boldsymbol{d}^{n+1} - \mathbf{g}(\boldsymbol{d}^n) + \dfrac{1}{\eta}(\nabla\phi^{n+1} - \boldsymbol{d}^{n+1}) + \tau(\boldsymbol{d}^n\cdot\boldsymbol{h})\cdot\boldsymbol{h}\Big), \\ \dfrac{\partial\boldsymbol{d}^{n+1}}{\partial n}|_{\partial\Omega} = 0, \end{cases} \tag{3.12}$$

with

$$\dot{\boldsymbol{d}}^{n+1} = \frac{\boldsymbol{d}^{n+1}-\boldsymbol{d}^n}{\delta t} + (\boldsymbol{u}_{\star\star}^n\cdot\nabla)\boldsymbol{d}^n,\ \boldsymbol{u}_{\star\star}^n = \boldsymbol{u}_\star^n - \delta t\frac{\dot{\boldsymbol{d}}^{n+1}}{M_2}\nabla\boldsymbol{d}^n. \tag{3.13}$$

*Step 3:*

$$\begin{cases} \dfrac{\tilde{\boldsymbol{u}}^{n+1} - \boldsymbol{u}^n}{\delta t} + (\boldsymbol{u}^n\cdot\nabla)\tilde{\boldsymbol{u}}^{n+1} - \mu_4\Delta\tilde{\boldsymbol{u}}^{n+1} + \nabla p^n + \dfrac{\dot{\phi}^{n+1}}{M_1}\nabla\phi^n + \dfrac{\dot{\boldsymbol{d}}^{n+1}}{M_2}\nabla\boldsymbol{d}^n = 0, \\ \tilde{\boldsymbol{u}}^{n+1}|_{\partial\Omega} = 0, \end{cases} \tag{3.14}$$

*Step 4:*

$$
\begin{cases}
\dfrac{\boldsymbol{u}^{n+1}-\tilde{\boldsymbol{u}}^{n+1}}{\delta t}+\nabla(p^{n+1}-p^n)=0,\\
\nabla\cdot\boldsymbol{u}^{n+1}=0,\\
\boldsymbol{n}\cdot\boldsymbol{u}^{n+1}|_{\partial\Omega}=0.
\end{cases}
\tag{3.15}
$$

In the above, $S$ is a stabilizing parameter to be determined.

We have the following remarks in order:

- A pressure-correction scheme [13] is used to decouple the computation of the pressure from that of the velocity.
- The nonlinear term $\mathbf{g}(\boldsymbol{d})$ mainly takes the form like $\frac{1}{\epsilon^2}\boldsymbol{d}(|\boldsymbol{d}|^2-1)$, so the explicit treatment of this term usually leads to a severe restriction on the time step $\delta t$ when $\epsilon \ll 1$. Thus we introduce in (3.10) a linear *"stabilizing"* term to improve the stability while preserving the simplicity. It allows us to treat the nonlinear term explicitly without suffering from any time step constraint [5,23,25,28–30,32–36,38, 40,41,45–47]. Note that this stabilizing term introduces an extra consistent error of order $O(\delta t)$ in a small region near the interface, but this error is of the same order as the error introduced by treating it explicitly, so the overall truncation error is essentially of the same order with or without the stabilizing term. It is remarkable that the truncation error of the stabilizing approach is exactly the same as the convex splitting method [9].
- Inspired by [1, 26, 31], which deal with a phase-field model of three-phase viscous fluids or complex fluids, we introduce two new, *explicit*, convective velocities $\boldsymbol{u}_\star^n$ and $\boldsymbol{u}_{\star\star}^n$ in the phase equations. $\boldsymbol{u}_\star^n$ and $\boldsymbol{u}_{\star\star}^n$ can be computed directly from (3.11) and (3.13), i.e.,

$$
\boldsymbol{u}_\star^n = \Big(I+\frac{\delta t}{M_1}(\nabla\phi^n)^T\nabla\phi^n\Big)^{-1}\Big(\boldsymbol{u}^n-\frac{1}{M_1}(\phi^{n+1}-\phi^n)\nabla\phi^n\Big), \tag{3.16}
$$

$$
\boldsymbol{u}_{\star\star}^n = \Big(I+\frac{\delta t}{M_2}(\nabla\boldsymbol{d}^n)^T\nabla\boldsymbol{d}^n\Big)^{-1}\Big(\boldsymbol{u}_\star^n-\frac{1}{M_2}(\boldsymbol{d}^{n+1}-\boldsymbol{d}^n)\nabla\boldsymbol{d}^n\Big). \tag{3.17}
$$

  It is easy to get $det(I+c(\nabla\phi)^T\nabla\phi)=1+c\nabla\phi\cdot\nabla\phi$, thus the above matrix is invertible.
- The scheme (3.10)-(3.15) is a totally decoupled, linear scheme. Indeed, (3.10), (3.12) and (3.14) are respectively (decoupled) linear elliptic equations for $\phi^{n+1}$, $\boldsymbol{d}^{n+1}$ and $\tilde{\boldsymbol{u}}^{n+1}$, and (3.15) can be reformulated as a Poisson equation for $p^{n+1}-p^n$. Therefore, at each time step, one only needs to solve a sequence of decoupled elliptic equations that can be solved very efficiently.
- As we shall show below, the above scheme is unconditionally energy stable.

**Theorem 3.1.** *Under the condition* (3.3)*, and* $S \geq \frac{\lambda L}{2}$*, the scheme* (3.10)-(3.15) *admits a unique solution satisfying the following discrete energy dissipation law:*

$$E^{n+1} + \frac{\delta t^2}{2}\|\nabla p^{n+1}\|^2 + \left\{\nu\delta t\|\nabla\tilde{\boldsymbol{u}}^{n+1}\|^2 + \delta t(\frac{|\dot{\phi}^{n+1}|^2}{M_1} + \frac{|\dot{\boldsymbol{d}}^{n+1}|^2}{M_2})\right\} \leq E^n + \frac{\delta t^2}{2}\|\nabla p^n\|^2,$$

*where*

$$E^n = \frac{1}{2}\|\boldsymbol{u}^n\|^2 + \lambda\Big(\eta\frac{\|\nabla\boldsymbol{d}^n\|^2}{2} + (\mathbf{G}(\boldsymbol{d}^n),1) + \frac{1}{\eta}\frac{\|\boldsymbol{d}^n - \nabla\phi^n\|^2}{2} - \frac{\tau}{2}\|\boldsymbol{d}^n\cdot\boldsymbol{h}\|^2\Big). \tag{3.18}$$

*Proof.* From the definition of $\boldsymbol{u}^n_\star$ and $\boldsymbol{u}^n_{\star\star}$ in (3.11) and (3.13), we can rewrite the momentum equation (3.14) as follows

$$\frac{\tilde{\boldsymbol{u}}^{n+1} - \boldsymbol{u}^n_{\star\star}}{\delta t} + (\boldsymbol{u}^n\cdot\nabla)\tilde{\boldsymbol{u}}^{n+1} - \mu_4\Delta\tilde{\boldsymbol{u}}^{n+1} + \nabla p^n = 0. \tag{3.19}$$

By taking the inner product of (3.19) with $2\delta t\tilde{\boldsymbol{u}}^{n+1}$, and using the identity

$$(a-b, 2a) = |a|^2 - |b|^2 + |a-b|^2, \tag{3.20}$$

we obtain

$$\|\tilde{\boldsymbol{u}}^{n+1}\|^2 - \|\boldsymbol{u}^n_{\star\star}\|^2 + \|\tilde{\boldsymbol{u}}^{n+1} - \boldsymbol{u}^n_{\star\star}\|^2 + 2\mu_4\delta t\|\nabla\tilde{\boldsymbol{u}}^{n+1}\|^2 + 2\delta t(\nabla p^n, \tilde{\boldsymbol{u}}^{n+1}) = 0. \tag{3.21}$$

To deal with the pressure term, we take the inner product of (3.15) with $2\delta t^2\nabla p^n$ to derive

$$\delta t^2(\|\nabla p^{n+1}\|^2 - \|\nabla p^n\|^2 - \|\nabla p^{n+1} - \nabla p^n\|^2) = 2\delta t(\tilde{\boldsymbol{u}}^{n+1}, \nabla p^n). \tag{3.22}$$

By taking the inner product of (3.15) with $\boldsymbol{u}^{n+1}$, we obtain

$$\|\boldsymbol{u}^{n+1}\|^2 + \|\boldsymbol{u}^{n+1} - \tilde{\boldsymbol{u}}^{n+1}\|^2 = \|\tilde{\boldsymbol{u}}^{n+1}\|^2. \tag{3.23}$$

We also derive from (3.15) directly that

$$\delta t^2\|\nabla p^{n+1} - \nabla p^n\|^2 = \|\tilde{\boldsymbol{u}}^{n+1} - \boldsymbol{u}^{n+1}\|^2. \tag{3.24}$$

Combining all identities above, we obtain

$$\begin{aligned}&\|\boldsymbol{u}^{n+1}\|^2 - \|\boldsymbol{u}^n_{\star\star}\|^2 + \|\tilde{\boldsymbol{u}}^{n+1} - \boldsymbol{u}^n_{\star\star}\|^2\\ &\qquad + \delta t^2(\|\nabla p^{n+1}\|^2 - \|\nabla p^n\|^2) + 2\nu\delta t\|\nabla\tilde{\boldsymbol{u}}^{n+1}\|^2 = 0.\end{aligned} \tag{3.25}$$

Next, we derive from (3.11) and (3.13) that

$$\frac{\boldsymbol{u}^n_\star - \boldsymbol{u}^n}{\delta t} = -\frac{\dot{\phi}^{n+1}}{M_1}\nabla\phi^n, \tag{3.26}$$

$$\frac{\boldsymbol{u}^n_{\star\star} - \boldsymbol{u}^n_\star}{\delta t} = -\frac{\dot{\boldsymbol{d}}^{n+1}}{M_2}\nabla\boldsymbol{d}^n. \tag{3.27}$$

By taking the inner product of (3.26) with $2\delta t\boldsymbol{u}_\star^n$, of (3.27) with $2\delta t\boldsymbol{u}_{\star\star}^n$, we obtain

$$\|\boldsymbol{u}_\star^n\|^2 - \|\boldsymbol{u}^n\|^2 + \|\boldsymbol{u}_\star^n - \boldsymbol{u}_n\|^2 = -2\delta t(\frac{\dot{\phi}^{n+1}}{M_1}\nabla\phi^n, \boldsymbol{u}_\star^n), \tag{3.28}$$

$$\|\boldsymbol{u}_{\star\star}^n\|^2 - \|\boldsymbol{u}_\star^n\|^2 + \|\boldsymbol{u}_{\star\star}^n - \boldsymbol{u}_\star^n\|^2 = -2\delta t(\frac{\dot{\boldsymbol{d}}^{n+1}}{M_2}\nabla\boldsymbol{d}^n, \boldsymbol{u}_{\star\star}^n). \tag{3.29}$$

Then, by taking the inner product of (3.10) with $2(\phi^{n+1} - \phi^n)$, we obtain

$$\begin{aligned} 2\delta t\frac{\|\dot{\phi}^{n+1}\|^2}{M_1} - 2\delta t(\frac{\dot{\phi}^{n+1}}{M_1}, (\boldsymbol{u}_\star^n\cdot\nabla)\phi^n) + \frac{\lambda}{\eta}\Big(\|\nabla\phi^{n+1}\|^2 - \|\nabla\phi^n\|^2 + \|\nabla\phi^{n+1} - \nabla\phi^n\|^2\Big) \\ + \frac{2\lambda}{\eta}\Big(\nabla\cdot\boldsymbol{d}^n, \phi^{n+1} - \phi^n\Big) = 0. \end{aligned} \tag{3.30}$$

By taking the inner product of (3.12) with $2(\boldsymbol{d}^{n+1} - \boldsymbol{d}^n)$, we arrive at

$$\begin{aligned} & 2S\|\boldsymbol{d}^{n+1} - \boldsymbol{d}^n\|^2 + 2\delta t\frac{\|\dot{\boldsymbol{d}}^{n+1}\|^2}{M_2} - 2\delta t(\frac{\dot{\boldsymbol{d}}^{n+1}}{M_2}, (\boldsymbol{u}_{\star\star}^n\cdot\nabla)\boldsymbol{d}^n) \\ & + 2\lambda\eta\Big(\frac{\|\nabla\boldsymbol{d}^{n+1}\|^2}{2} - \frac{\|\nabla\boldsymbol{d}^n\|^2}{2} + \frac{\|\nabla\boldsymbol{d}^{n+1} - \nabla\boldsymbol{d}^n\|^2}{2}\Big) \\ & + \frac{\lambda}{\eta}(\|\boldsymbol{d}^{n+1}\|^2 - \|\boldsymbol{d}^n\|^2 + \|\boldsymbol{d}^{n+1} - \boldsymbol{d}^n\|^2) \\ & + 2\lambda(\mathbf{g}(\boldsymbol{d}^n), \boldsymbol{d}^{n+1} - \boldsymbol{d}^n) - \frac{2\lambda}{\eta}\Big(\nabla\phi^{n+1}, \boldsymbol{d}^{n+1} - \boldsymbol{d}^n\Big) \\ & - 2\lambda\tau\Big((\boldsymbol{d}^n\cdot\boldsymbol{h})\boldsymbol{h}, \boldsymbol{d}^{n+1} - \boldsymbol{d}^n\Big) = 0. \end{aligned} \tag{3.31}$$

Combining (3.25), (3.28), (3.29), (3.30), and (3.31), we arrive at

$$
\begin{aligned}
&\|\boldsymbol{u}^{n+1}\|^2 - \|\boldsymbol{u}^n\|^2 + \|\tilde{\boldsymbol{u}}^{n+1} - \boldsymbol{u}^n_{\star\star}\|^2 + \|\boldsymbol{u}^n_{\star\star} - \boldsymbol{u}^n_{\star}\|^2 + \|\boldsymbol{u}^n_{\star} - \boldsymbol{u}^n\|^2 \\
&\quad + \delta t^2(\|\nabla p^{n+1}\|^2 - \|\nabla p^n\|^2) \\
&\quad + 2\nu\delta t\|\nabla\tilde{\boldsymbol{u}}^{n+1}\|^2 + 2\delta t\frac{\|\dot{\phi}^{n+1}\|^2}{M_1} + 2\delta t\frac{\|\dot{\boldsymbol{d}}^{n+1}\|^2}{M_2} \\
&\quad + 2\lambda\eta\Big(\frac{\|\nabla\boldsymbol{d}^{n+1}\|^2}{2} - \frac{\|\nabla\boldsymbol{d}^n\|^2}{2} + \frac{\|\nabla\boldsymbol{d}^{n+1} - \nabla\boldsymbol{d}^n\|^2}{2}\Big) \\
&\quad + \frac{\lambda}{\eta}(\|\boldsymbol{d}^{n+1}\|^2 - \|\boldsymbol{d}^n\|^2 + \|\boldsymbol{d}^{n+1} - \boldsymbol{d}^n\|^2) \\
&\quad + \frac{\lambda}{\eta}(\|\nabla\phi^{n+1}\|^2 - \|\nabla\phi^n\|^2 + \|\nabla\phi^{n+1} - \nabla\phi^n\|^2) \\
&\quad + 2S\|\boldsymbol{d}^{n+1} - \boldsymbol{d}^n\|^2 \\
&\quad + 2\lambda(\mathbf{g}(\boldsymbol{d}^n), \boldsymbol{d}^{n+1} - \boldsymbol{d}^n) \text{ (:Term } A) \\
&\quad + \frac{2\lambda}{\eta}(\nabla\cdot\boldsymbol{d}^n, \phi^{n+1} - \phi^n) - \frac{2\lambda}{\eta}(\nabla\phi^{n+1}, \boldsymbol{d}^{n+1} - \boldsymbol{d}^n) \text{ (:Term } B) \\
&\quad - 2\lambda\tau\Big((\boldsymbol{d}^n\cdot\boldsymbol{h})\boldsymbol{h}, \boldsymbol{d}^{n+1} - \boldsymbol{d}^n\Big) \text{ (:Term } C) \\
&\quad = 0.
\end{aligned}
\tag{3.32}
$$

We deal with the terms $A, B, C$ as follows.

For Term $A$, we apply the Taylor expansions to obtain

$$
A = 2\lambda(\mathbf{G}(\boldsymbol{d}^{n+1}) - \mathbf{G}(\boldsymbol{d}^n), 1) - 2\lambda\Big(\frac{\mathbf{g}_1'(\xi)}{2}, |\boldsymbol{d}^{n+1} - \boldsymbol{d}^n|^2\Big). \tag{3.33}
$$

For Term $B$, we have

$$
\begin{aligned}
B &= -2\frac{\lambda}{\eta}\Big((\boldsymbol{d}^n, \nabla\phi^{n+1} - \nabla\phi^n) + (\boldsymbol{d}^{n+1} - \boldsymbol{d}^n, \nabla\phi^{n+1})\Big) \\
&= -2\frac{\lambda}{\eta}\Big((\boldsymbol{d}^{n+1}, \nabla\phi^{n+1}) - (\boldsymbol{d}^n, \nabla\phi^n)\Big).
\end{aligned}
\tag{3.34}
$$

For Term $C$, we have

$$
\begin{aligned}
C &= -2\lambda\tau\Big((\boldsymbol{d}^n\cdot\boldsymbol{h}), (\boldsymbol{d}^{n+1}\cdot\boldsymbol{h}) - (\boldsymbol{d}^n\cdot\boldsymbol{h})\Big) \\
&= -\lambda\tau\Big(\|\boldsymbol{d}^{n+1}\cdot\boldsymbol{h}\|^2 - \|\boldsymbol{d}^n\cdot\boldsymbol{h}\|^2 - \|(\boldsymbol{d}^{n+1} - \boldsymbol{d}^n)\cdot\boldsymbol{h}\|^2\Big).
\end{aligned}
\tag{3.35}
$$

By combining (3.32), (3.33), (3.34) and (3.35), we have

$$
\begin{aligned}
&\|\boldsymbol{u}^{n+1}\|^2 - \|\boldsymbol{u}^n\|^2 + \|\tilde{\boldsymbol{u}}^{n+1} - \boldsymbol{u}^n_{\star\star}\|^2 + \|\boldsymbol{u}^n_{\star\star} - \boldsymbol{u}^n_{\star}\|^2 + \|\boldsymbol{u}^n_{\star} - \boldsymbol{u}^n\|^2 \\
&\quad + \delta t^2(\|\nabla p^{n+1}\|^2 - \|\nabla p^n\|^2) \\
&\quad + 2\nu\delta t\|\nabla \tilde{\boldsymbol{u}}^{n+1}\|^2 + 2\delta t\frac{\|\dot{\phi}^{n+1}\|^2}{M_1} + 2\delta t\frac{\|\dot{\boldsymbol{d}}^{n+1}\|^2}{M_2} \\
&\quad + 2\lambda\eta\Big(\frac{\|\nabla \boldsymbol{d}^{n+1}\|^2}{2} - \frac{\|\nabla \boldsymbol{d}^{n}\|^2}{2} + \frac{\|\nabla \boldsymbol{d}^{n+1} - \nabla \boldsymbol{d}^n\|^2}{2}\Big) \\
&\quad + 2\lambda(\mathbf{G}(\boldsymbol{d}^{n+1}) - \mathbf{G}(\boldsymbol{d}^n), 1) \\
&\quad + \frac{\lambda}{\eta}(\|\boldsymbol{d}^{n+1}\|^2 - \|\boldsymbol{d}^n\|^2 + \|\boldsymbol{d}^{n+1} - \boldsymbol{d}^n\|^2) \\
&\quad + \frac{\lambda}{\eta}(\|\nabla\phi^{n+1}\|^2 - \|\nabla\phi^n\|^2 + \|\nabla\phi^{n+1} - \nabla\phi^n\|^2) \\
&\quad - 2\frac{\lambda}{\eta}\Big((\boldsymbol{d}^{n+1}, \nabla\phi^{n+1}) - (\boldsymbol{d}^n, \nabla\phi^n)\Big) \\
&\quad - \lambda\tau\Big(\|\boldsymbol{d}^{n+1}\cdot \boldsymbol{h}\|^2 - \|\boldsymbol{d}^n \cdot \boldsymbol{h}\|^2\Big) \\
&\quad + (2S - \lambda L)\|\boldsymbol{d}^{n+1} - \boldsymbol{d}^n\|^2 \\
&\quad + \lambda\tau\|(\boldsymbol{d}^{n+1} - \boldsymbol{d}^n)\cdot \boldsymbol{h}\|^2 = 0.
\end{aligned} \tag{3.36}
$$

Finally, we obtain the desired result after dropping some positive terms. □

## 4. Numerical Simulations

We now present several 2D numerical experiments to demonstrate the efficiency, stability and accuracy of the propose numerical scheme (3.10)-(3.15). The computational domain is $(x, y) \in [0, 2\ell] \times [0, 2]$. For $x$-axis, we set the periodic boundary condition, and for $y-$ axis, we set Neumann or Dirichlet boundary conditions. We adopt the spectral Galerkin method to discretize the space. For $x$-axis, we use the Fourier basis, and the for $y - axis$, we use the Legendre polynomial (cf. [27]). The magnetic field is always set as $\boldsymbol{h} = (0, 1)$. If not explicitly specified, the default values of order parameters are given as follows,

$$\ell = 2, \epsilon = 0.02,\ \eta = 0.02,\ M_1 = 0.08,\ M_2 = 2,\ \lambda = 2.5,\ \tau = 16. \tag{4.1}$$

4.1. **Accuracy test.** We first test the convergence rate of scheme (3.10)-(3.15). We set the following initial conditions as

$$
\begin{cases}
\boldsymbol{d}(t=0) = (\sin(\pi x)\cos(\pi y), \cos(\pi x)\cos(\pi y)), \\
\phi(t=0) = \cos(\pi y), \\
(\boldsymbol{u}(t=0), p(t=0)) = (0, 0).
\end{cases} \tag{4.2}
$$

The boundary conditions are Neumann type along y-axis (cf. (2.15)). We use $128 \times 128$ grid points to discretize the space, and perform the mesh refinement test for time. We choose the numerical solution with the time step size $\delta t = 1 \times 10^{-4}$ as the benchmark solution (approximate exact solution) for computing errors. Figure 4.1 plots the $L^2$ errors for various time step sizes. We observe that the scheme is asymptotically first-order accurate in time for all variables as expected.

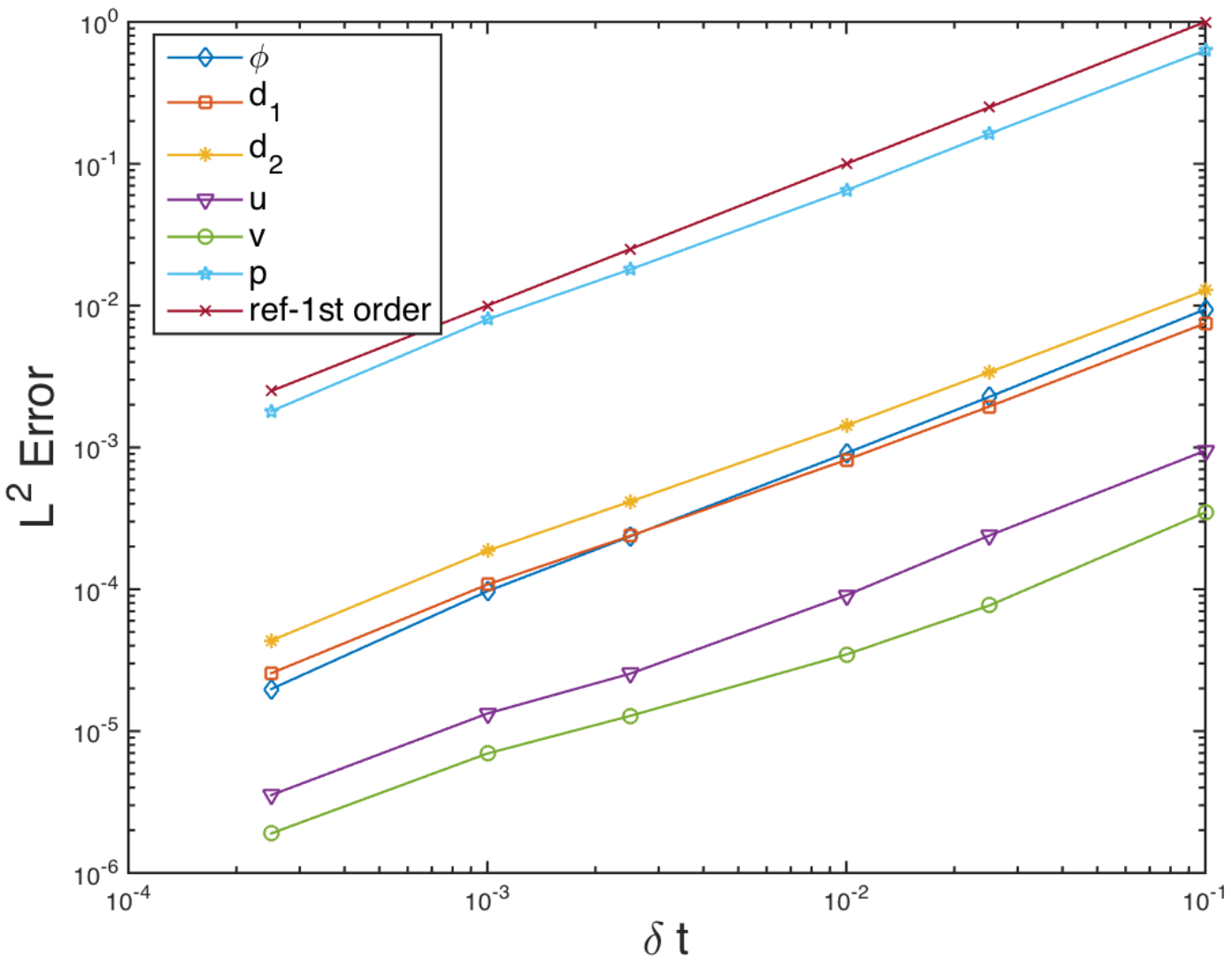

FIGURE 4.1. The $L^2$ errors of the layer funciton $\phi$, the director field $\boldsymbol{d} = (d_1, d_2)$, the velocity $\boldsymbol{u} = (u, v)$ and pressure $p$. The slopes show that the scheme is asymptotically first-order accurate in time.

4.2. **Chevron pattern induced by the magnetic force.** We now consider the effects from the magnetic force for the no flow case ($\boldsymbol{u} = 0$). Initially, a smectic A liquid crystal is confined between two flat parallel plates and uniformly aligned in a way that the smectic layers are parallel to the bounding plates and the directors are aligned homeotropically, that is, perpendicular to the smectic layers. A magnetic field is applied in the direction parallel to the smectic layers, which induce the layer undulation (chevron pattern) phenomena. The initial conditions read as follows.

$$\begin{aligned} &\boldsymbol{d}(t=0) = (0,1) + 0.001(\mathrm{rand}(x,y), \mathrm{rand}(x,y)), \\ &\phi(t=0) = y, \end{aligned} \tag{4.3}$$

where the $\mathrm{rand}(x, y)$ is the small perturbation that is the random number in $[-1, 1]$ and has zero mean. We set the Dirichlet type boundary condition for $\phi$ and $\boldsymbol{d}$ as follows,

$$\boldsymbol{d}|_{y=\pm1} = (0,1),\ \phi|_{y=1} = 1,\ \phi|_{y=-1} = -1. \tag{4.4}$$

We take $\delta t = 0.001$ to obtain better accuracy. Fig. 4.2 shows the snapshots of the layer function $\phi$ at $t = 0$, 0.2, 0.4 and 0.8. Initially at $t = 0$, the layer function take the linear profile along the $y-$ axis. When time evolves, we observe some undulations appear at $t = 0.2$. The layer function quickly reaches the steady solution at $t = 0.8$ with the saw tooth shape. This undulation phenomenon is called the Helfrich-Hurault effect (cf. [15,16]). Fig. 4.3 shows the snapshots of the director field $\boldsymbol{d}$. The numerical solution presents similar features to those obtained in [18]. We also plot the energy dissipative curve in Fig. 4.4, which confirms that our algorithm is energy stable.

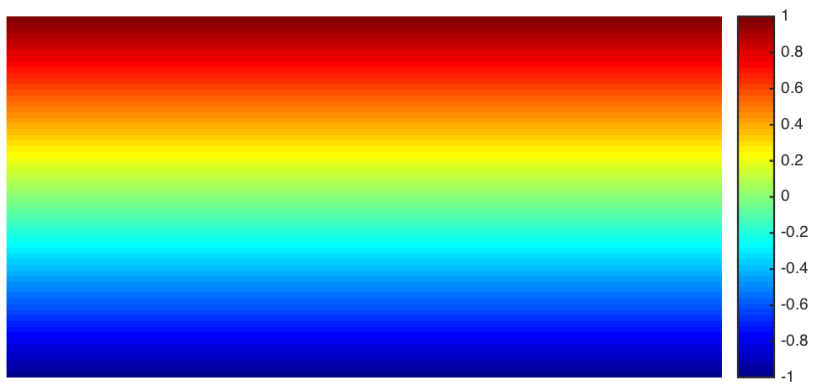

(a) $t = 0$.

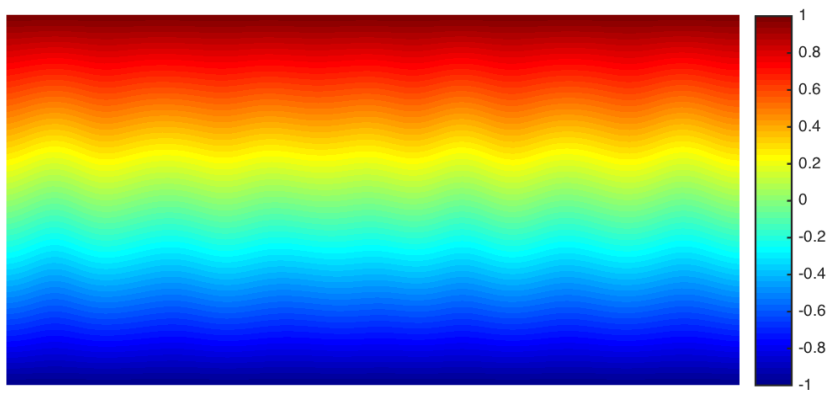

(b) $t = 0.2$.

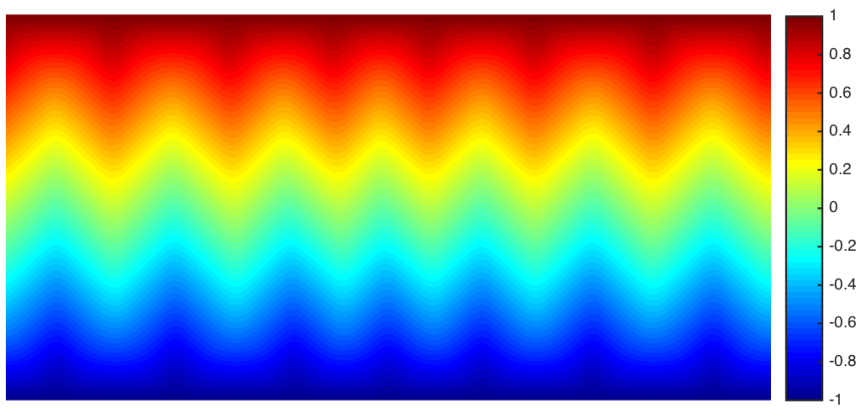

(c) $t = 0.4$.

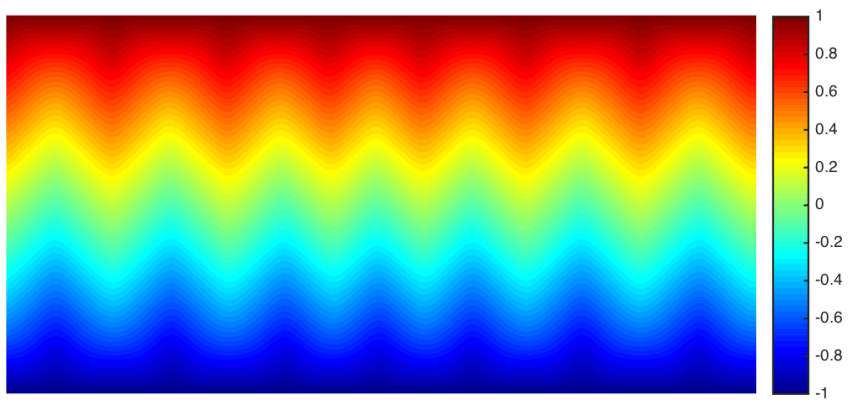

(d) $t = 0.8$.

FIGURE 4.2. Snapshots of the layer function $\phi$ are taken at $t = 0$, 0.2, 0.4 and 0.8 for Example 4.2.

4.3. **Chevron pattern induced by magnetic force and shear flow.** We now impose the shear flow on the top and bottom plates to see how the flow affects the undulation. The initial and boundary conditions of $\phi$ and $\boldsymbol{d}$ are same as the example 4.2. For velocity and pressure, the initial and boundary conditions are:

$$\begin{aligned} &\boldsymbol{u}(t=0) = (10(y-1), 0),\ p(t=0) = 0 \\ &\boldsymbol{u}|_{y=1} = (10, 0), \boldsymbol{u}|_{y=-1} = (-10, 0). \end{aligned} \tag{4.5}$$

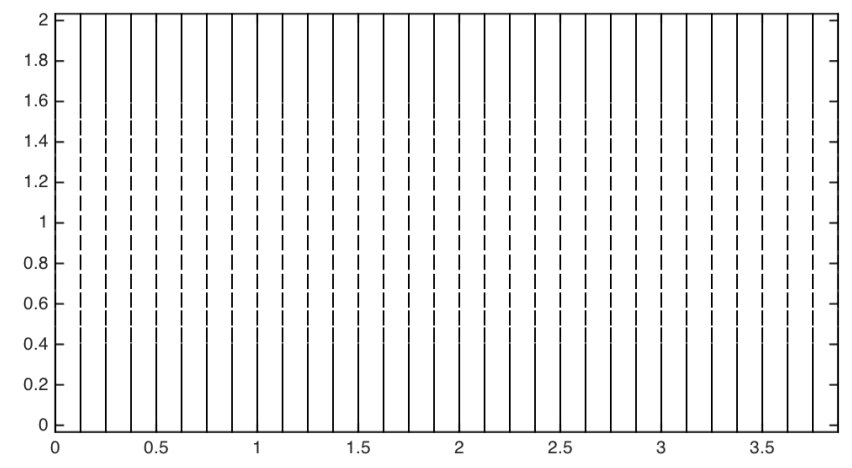

(a) $t = 0$.

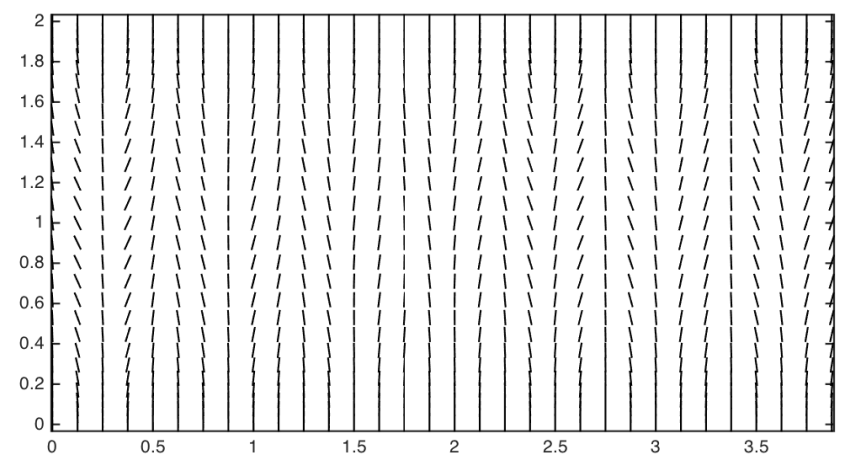

(b) $t = 0.2$.

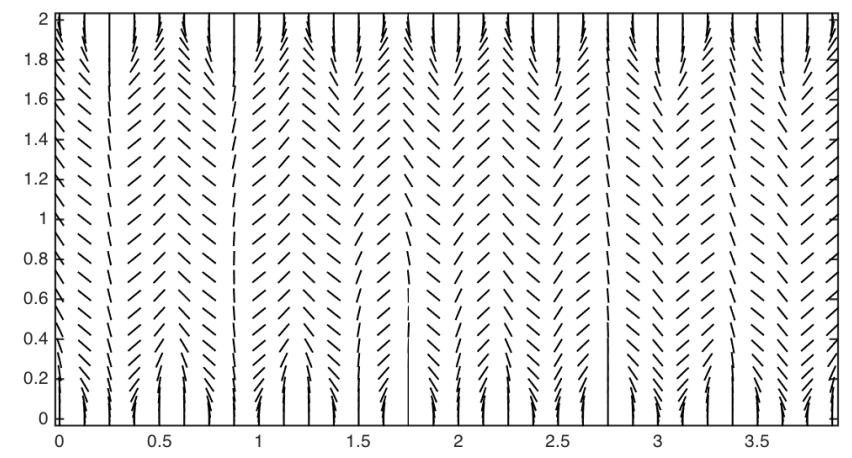

(c) $t = 0.4$.

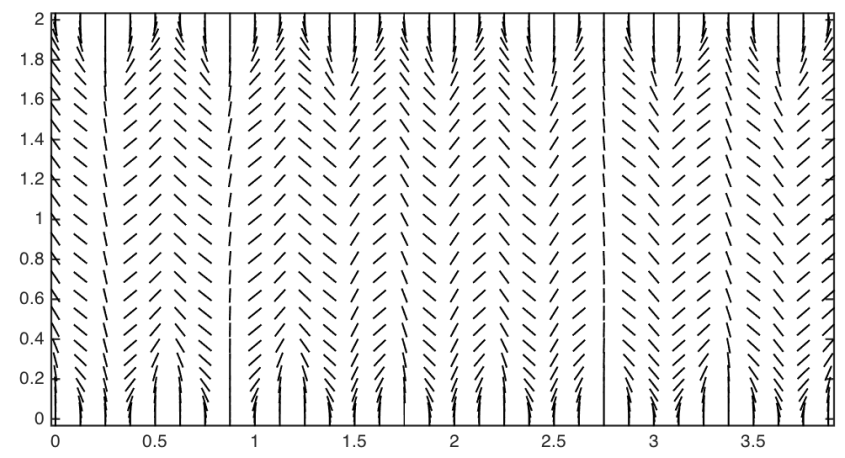

(d) $t = 0.8$.

FIGURE 4.3. Snapshots of the director field $\boldsymbol{d}$ are taken at $t = 0$, $0.2$, $0.4$ and $0.8$ for Example 4.2.

Fig. 4.5 show the snapshots of the layer function $\phi$ at $t = 0, 0.3, 0.4, 0.5, 0.6$ and $0.8$. When time evolves, the layer undulations still appear but the symmetry is largely disturbed by the shear flow. Fig. 4.6 shows the snapshots of the directior field $\boldsymbol{d}$. We also plot the first component of the velocity field $\boldsymbol{u} = (u, v)$ in Fig. 4.7, where the initial linear profile is deformed and the nonlinearity is shown.

## 5. Concluding Remarks

In this paper, we have presented an efficient, semi-discrete in time, numerical scheme with provably unconditionally stability for solving the hydrodynamics coupled liquid crystal flows of smectic-A phase. Based on the stabilized approache, the proposed scheme conquer the inconvenience from nonlinearities by linearizing all nonlinear terms explicitly. We show that the linear scheme developed is unconditionally energy stable. We verify numerically that our scheme is of first order accurate in time and present ample numerical results from some standard numerical tests in the presence of the magnetic field and the shear flow.

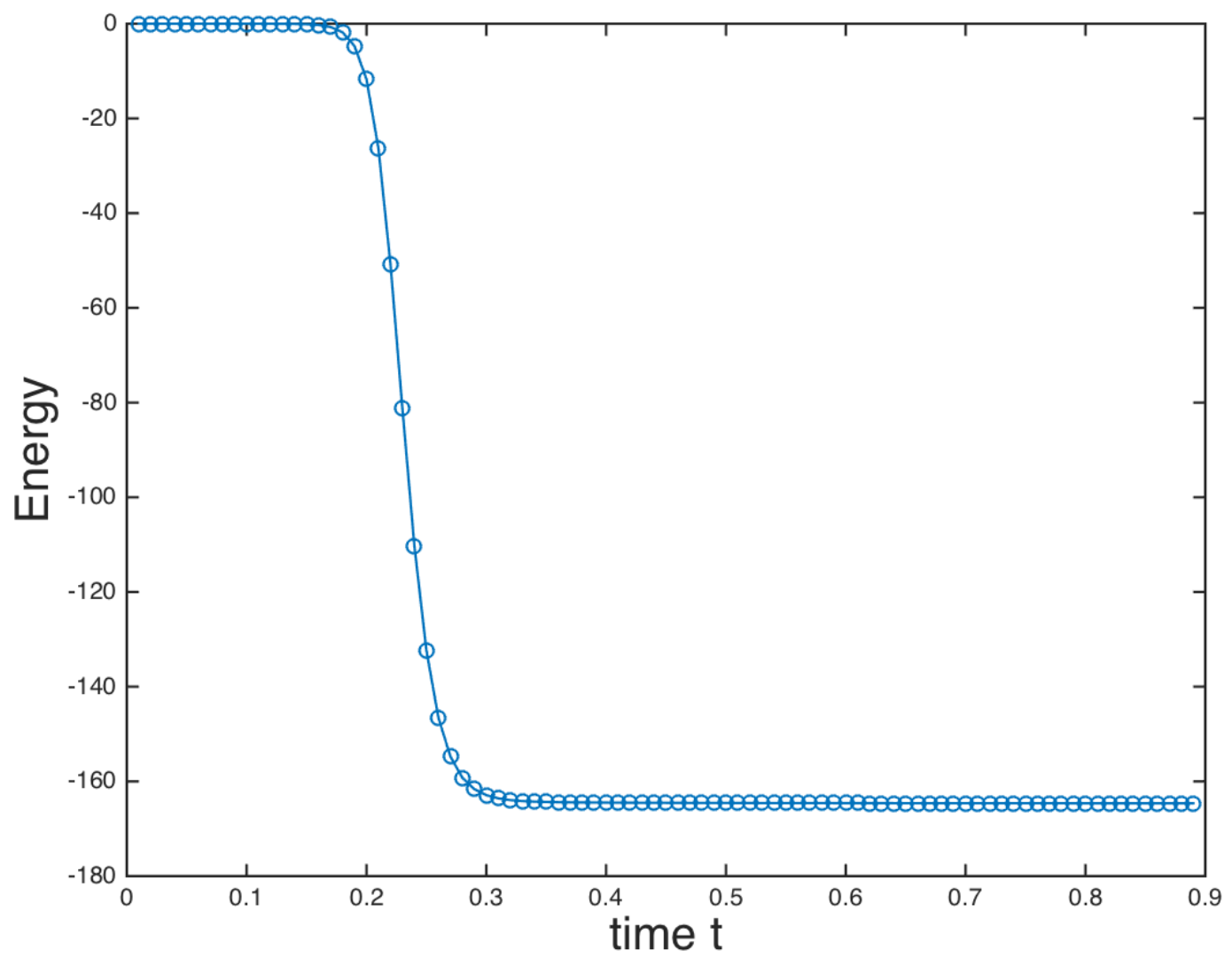

Figure 4.4. Time evolution of the free energy functional of Example 4.2.

**Acknowledgments.** X. Yang's research is partially supported by the U.S. National Science Foundation under grant numbers DMS-1418898.

## References

[1] F. Boyer and S. Minjeaud. Numerical schemes for a three component Cahn-Hilliard model. *ESAIM Math. Model. Numer. Anal.*, 45(4):697–738, 2011.

[2] J. W. Cahn and J. E. Hilliard. Free energy of a nonuniform system. I. interfacial free energy. *J. Chem. Phys.*, 28:258–267, 1958.

[3] S. Chandrasekhar. Liquid crystals (2nd ed.). cambridge: Cambridge university press, isbn 0-521-41747-3. 1992.

[4] J. Chen and T. C. Lubensky. Landau-ginzburg mean-field theory for the nematic to smectic-c and nematic to smectic-a phase transitions. *Phys. Rev. A.*, 14:1202–1207, 1976.

[5] R. Chen, G. Ji, X. Yang, and H. Zhang. Decoupled energy stable schemes for phase-field vesicle membrane model. *J. Comput. Phys.*, 302:509–523, 2015.

[6] N. Condette, C. Melcher, and E. Süli. Spectral approximation of pattern-forming nonlinear evolution equations with double-well potentials of quadratic growth. *to appear in Math. Comp.*

[7] P. G. de Gennes and J. Prost. *The Physics of Liquid Crystals.* Oxford University Press, 1993.

[8] W. E. Nonlinear continuum theory of smectic-a liquid crystals. *Arch. Ration. Mech. Anal.*, 137:159–175, 1997.

[9] D. J. Eyre. Unconditionally gradient stable time marching the Cahn-Hilliard equation. In *Computational and mathematical models of microstructural evolution (San Francisco, CA, 1998)*, volume 529 of *Mater. Res. Soc. Sympos. Proc.*, pages 39–46. MRS, Warrendale, PA, 1998.

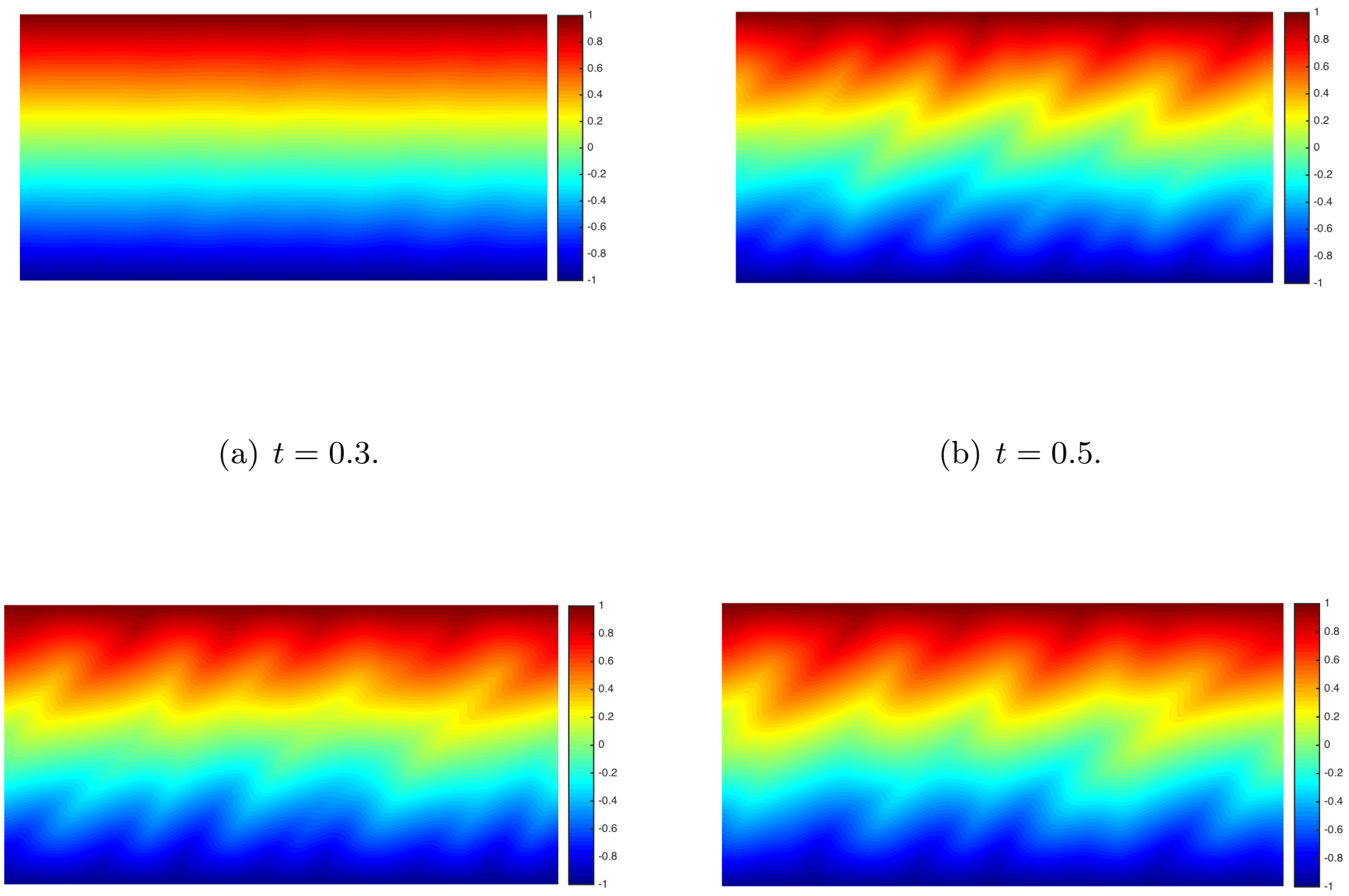

(a) $t = 0.3$.

(b) $t = 0.5$.

(c) $t = 0.6$.

(d) $t = 0.8$.

FIGURE 4.5. Snapshots of the layer function $\phi$ are taken at $t = 0.3$ , 0.5, 0.6 and 0.8 for Example 4.3.

[10] A. Fick. Poggendorff's annalen. *Journal of the american mathematics society*, pages 59–86, 1855.

[11] M. G. Forest, S. Heidenreich, S. Hess, X. Yang, and R. Zhou. Robustness of pulsating jetlike layers in sheared nano-rod dispersions. *J. Non-Newtonian Fluid Mech.*, 155:130–145, 2008.

[12] M. G. Forest, S. Heidenreich, S. Hess, X. Yang, and R. Zhou. Dynamic texture scaling of sheared nematic polymers in the large ericksen number limit. *J. Non-Newtonian Fluid Mech.*, 165:687–697, 2010.

[13] J. L. Guermond, P. Minev, and J. Shen. An overview of projection methods for incompressible flows. *Comput. Methods Appl. Mech. Engrg.*, 195:6011–6045, 2006.

[14] F. Guillen-Gonzaleza and G. Tierra. Approximation of smectic-a liquid crystals. *Comput. Methods Appl. Mech. Engrg.*, 290:342–361, 2015.

[15] W. Helfrich. Electrohydrodynamic and dielectric instabilities of cholesteric liquid crystals. *The Journal of Chemical Physics*, 55:839–842, 1971.

[16] J. P. Hurault. Static distortions of a cholesteric planar structure induced by magnetic or ac electric fields. *The Journal of Chemical Physics*, 59:2068–2075, 1973.

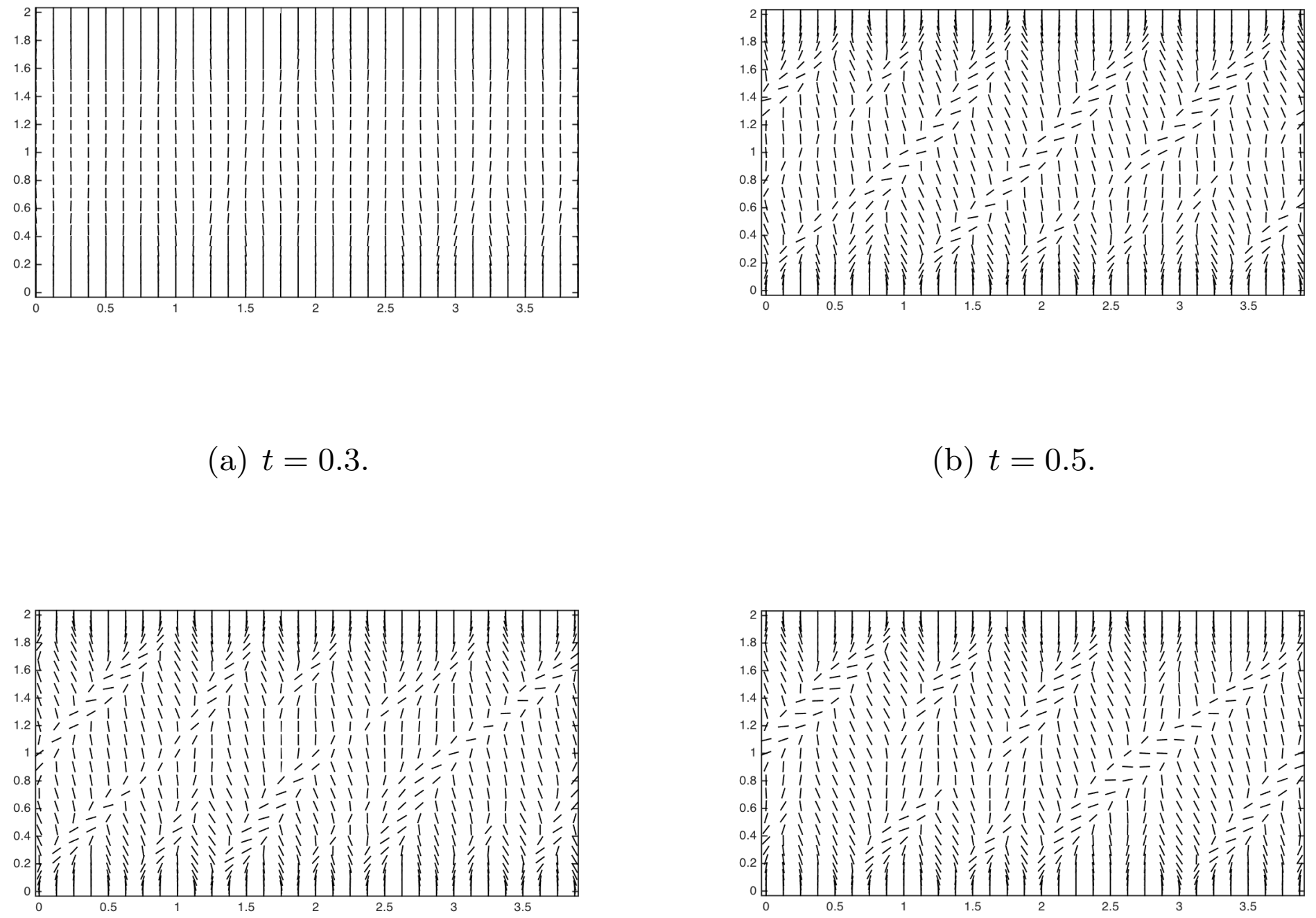

(a) $t = 0.3$.

(b) $t = 0.5$.

(c) $t = 0.6$.

(d) $t = 0.8$.

FIGURE 4.6. Snapshots of the director field $\boldsymbol{d}$ are taken at $t = 0.3$, 0.5, 0.6 and 0.8 for Example 4.3.

[17] Y. Gong J. Zhao, X. Yang and Q. Wang. A novel linear second order unconditionally energy stable scheme for a hydrodynamic q-tensor model of liquid crystals. *in press, Comput. Meth. Appl. Mech. Engrg.*, 2017.

[18] S. Joo and D. Phillips. The phase transitions from chiral nematic toward smectic liquid crystals. *Communications in Mathematical Physics*, 269:367–399, 2007.

[19] D. Kessler, R. H. Nochetto, and A. Schmidt. A posteriori error control for the Allen-Cahn problem: circumventing Gronwall's inequality. *M2AN Math. Model. Numer. Anal.*, 38(1):129–142, 2004.

[20] F. H. Lin. On nematic liquid crystals with variable degree of orientation. *Communications on Pure and Applied Mathematics*, 44:453–468, 1991.

[21] F. H. Lin. Mathematics theory of liquid crystals, in applied mathematics at the turn of century: Lecture notes of the 1993 summer school, universidat complutense de madrid. 1995.

[22] C. Liu and J. Shen. A phase field model for the mixture of two incompressible fluids and its approximation by a Fourier-spectral method. *Physica D*, 179(3-4):211–228, 2003.

[23] C. Liu, J. Shen, and X. Yang. Decoupled energy stable schemes for a phase-field model of two-phase incompressible flows with variable density. *J. Sci. Comput.*, 62:601–622, 2015.

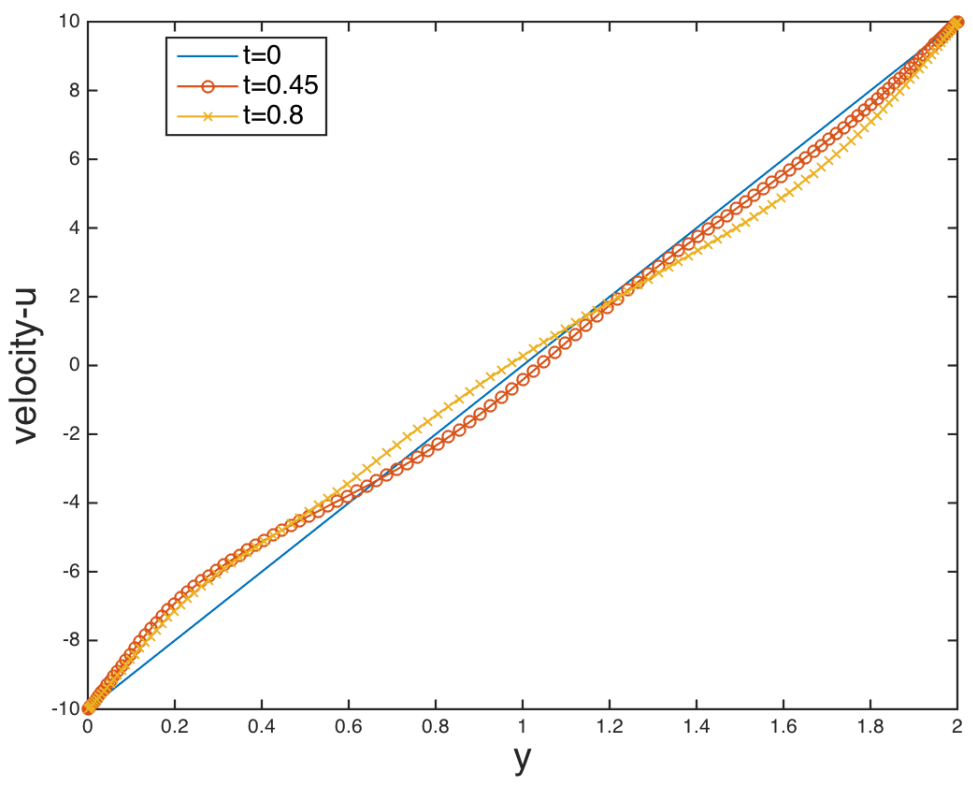

FIGURE 4.7. Snapshots of the profile for the first component $u(y)$ of the velocity field $\boldsymbol{u} = (u, v)$ at the center $(x = 2)$ and $t = 0, 0.45$ and 0.8.

[24] C. Liu and N.J. Walkington. An Eulerian description of fluids containing visco-hyperelastic particles. *Arch. Rat. Mech. Anal.*, 159:229–252, 2001.

[25] L. Ma, R. Chen, X. Yang, and H. Zhang. Numerical approximations for allen-cahn type phase field model of two-phase incompressible fluids with moving contact lines. *Comm. Comput. Phys.*, 21:867–889, 2017.

[26] S. Minjeaud. An unconditionally stable uncoupled scheme for a triphasic cahn-hilliard/navier-stokes model. *Commun. Comput. Phys.*, 29:584–618, 2013.

[27] J. Shen. Efficient spectral-Galerkin method I. direct solvers for second- and fourth-order equations by using Legendre polynomials. *SIAM J. Sci. Comput.*, 15:1489–1505, 1994.

[28] J. Shen and X. Yang. An efficient moving mesh spectral method for the phase-field model of two-phase flows. *J. Comput. Phys.*, 228:2978–2992, 2009.

[29] J. Shen and X. Yang. Energy stable schemes for Cahn-Hilliard phase-field model of two-phase incompressible flows. *Chin. Ann. Math. Ser. B*, 31(5):743–758, 2010.

[30] J. Shen and X. Yang. Numerical Approximations of Allen-Cahn and Cahn-Hilliard Equations. *Disc. Conti. Dyn. Sys.-A*, 28:1669–1691, 2010.

[31] J. Shen and X. Yang. A phase-field model and its numerical approximation for two-phase incompressible flows with different densities and viscositites. *SIAM J. Sci. Comput.*, 32:1159–1179, 2010.

[32] J. Shen and X. Yang. Decoupled energy stable schemes for phase filed models of two phase complex fluids. *SIAM J. Sci. Comput.*, 36:B122–B145, 2014.

[33] J. Shen and X. Yang. Decoupled, energy stable schemes for phase-field models of two-phase incompressible flows. *SIAM J. Num. Anal.*, 53(1):279–296, 2015.

[34] J. Shen, X. Yang, and Q. Wang. On mass conservation in phase field models for binary fluids. *Comm. Compt. Phys*, 13:1045–1065, 2012.

[35] J. Shen, X. Yang, and H. Yu. Efficient energy stable numerical schemes for a phase field moving contact line model. *J. Comput. Phys.*, 284:617–630, 2015.

[36] C. Xu and T. Tang. Stability analysis of large time-stepping methods for epitaxial growth models. *SIAM. J. Num. Anal.*, 44:1759–1779, 2006.

[37] K. Xu, M. G. Forest, and X. Yang. Shearing the i-n phase transition of liquid crystalline polymers: long-time memory of defect initial data. *Disc. Conti. Dyn. Sys.-B*, 15:457–474, 2010.

[38] X. Yang. Error analysis of stabilized semi-implicit method of Allen-Cahn equation. *Disc. Conti. Dyn. Sys.-B*, 11:1057–1070, 2009.

[39] X. Yang, Z. Cui, M. G. Forest, Q. Wang, and J. Shen. Dimensional robustness & instability of sheared, semi-dilute, nano-rod dispersions. *SIAM Multi. Model. Simul.*, 7:622–654, 2008.

[40] X. Yang, J. J. Feng, C. Liu, and J. Shen. Numerical simulations of jet pinching-off and drop formation using an energetic variational phase-field method. *J. Comput. Phys.*, 218:417–428, 2006.

[41] X. Yang, M. G. Forest, H. Li, C. Liu, J. Shen, Q. Wang, and F. Chen. Modeling and simulations of drop pinch-off from liquid crystal filaments and the leaky liquid crystal faucet immersed in viscous fluids. *J. Comput. Phys.*, 236:1–14, 2013.

[42] X. Yang, M. G. Forest, W. Mullins, and Q. Wang. 2-d lid-driven cavity flow of nematic polymers: an unsteady sea of defects. *Soft. Matter*, 6:1138–1156, 2009.

[43] X. Yang, M. G. Forest, W. Mullins, and Q. Wang. Dynamic defect morphology and hydrodynamics of sheared nematic polymers in two space dimensions. *J. Rheoloogy*, 53:589–615, 2009.

[44] X. Yang, M. G. Forest, W. Mullins, and Q. Wang. Quench sensitivity to defects and shear banding in nematic polymer film flows. *J. Non-Newtonian Fluid Mech.*, 159:115–129, 2009.

[45] H. Yu and X. Yang. Numerical approximations for a phase-field moving contact line model with variable densities and viscosities. *in press, DOI: 10.1016/j.jcp.2017.01.026, J. Comput. Phys.*, 2017.

[46] J. Zhao, H. Li, Q. Wang, and X. Yang. A linearly decoupled energy stable scheme for phase-field models of three-phase incompressible flows. *in press, DOI: 10.1007/s10915-016-0283-9, J. Sci. Comput.*, 2017.

[47] J. Zhao, Q. Wang, and X. Yang. Numerical approximations to a new phase field model for immiscible mixtures of nematic liquid crystals and viscous fluids. *Comput. Meth. Appl. Mech. Eng.*, 310:77–97, 2016.

[48] J. Zhao, X. Yang, J. Li, and Q. Wang. Energy stable numerical schemes for a hydrodynamic model of nematic liquid crystals. *SIAM. J. Sci. Comput.*, 38:A3264–A3290, 2016.

[49] J. Zhao, X. Yang, J. Shen, and Q. Wang. A decoupled energy stable scheme for a hydrodynamic phase-field model of mixtures of nematic liquid crystals and viscous fluids. *J. Comput. Phys.*, 305:539–556, 2016.